\documentclass[a4paper]{article}
\usepackage{geometry}
\usepackage{PRIMEarxiv}
\usepackage{microtype}
\usepackage[T1]{fontenc}
\usepackage[utf8]{inputenc}
\usepackage{amsmath,amssymb,amsthm}
\usepackage{algorithm}
\usepackage{algpseudocode}
\usepackage{float}
\usepackage{booktabs}
\usepackage{xcolor}
\usepackage{hyperref}
\usepackage{tikz}
\usetikzlibrary{arrows.meta,positioning,calc,shapes.geometric,fit,decorations.pathmorphing,decorations.pathreplacing}

\hypersetup{
  colorlinks=true,
  linkcolor=blue!50!black,
  citecolor=blue!50!black,
  urlcolor=blue!50!black
}

\newtheorem{remark}{Remark}
\newtheorem{proposition}{Proposition}

\title{\textbf{Optimizing Explicit Unit-Distance Lower-Bound Certificates}}
\author{Michael T. M. Emmerich\thanks{Faculty of Information Technology, University of Jyv\"askyl\"a, Finland. ORCID: \href{https://orcid.org/0000-0002-7342-2090}{0000-0002-7342-2090}.}}
\date{May 30, 2026}

\begin{document}
\maketitle

\begin{abstract}
The 2026 disproof of Erd\H{o}s's unit-distance conjecture and Sawin's quantitative refinement show that the maximum number $u(n)$ of unit distances among $n$ planar points can exceed $n^{1+\varepsilon}$ for a fixed positive $\varepsilon$. Sawin's explicit bound gives more than $n^{1.014}$ unit distances for arbitrarily large $n$ and exposes integer parameters whose choice is not fully optimized. This report treats Sawin's parameter selection as a nonlinear integer optimization problem and develops an open-source Python optimization and verification pipeline for certificates involving prime sets $T$ and $S_Q$, integer multiplicities $k(p)$, and a rationally encoded real parameter $R$. After reproducing Sawin's certificate with $\delta=0.014114\ldots$, the pipeline yields improved certificates with the same $T$. We develop a tailored integer evolution strategy achieving a certificate with $\delta=0.015263\ldots$ and supporting the cautious statement $u(n)>n^{1.0152}$ for arbitrarily large $n$. For extended ramified prime ranges, the Emmerich--Cordella certificate obtained with the same framework reports $u(n)>n^{1.031}$ for $\#T=67$, illustrating the importance of enlarging $T$. Very recent MathOverflow discussions, brought to the author's attention as of version~4, report further improvements, including certificates above $\delta>0.035$ and beyond $\delta>0.036$. Some of these improvements may rely not only on larger prime ranges but also on modified constraint systems and additional degrees of freedom that deviate from Sawin's original formulation. Beyond this application, the work illustrates how randomized optimization heuristics can improve, verify, and refine explicit certificates for combinatorial geometry through nonlinear integer optimization.
\end{abstract}

\keywords{Erd\H{o}s unit-distance problem; discrete geometry; algebraic number theory; Golod--Shafarevich inequality; class-field towers; nonlinear integer programming; integer evolution strategy; reproducible verification.}

\medskip

\noindent\textbf{MSC 2020.} 52C10; 11R29; 11R37; 90C11; 90C59.

\medskip

\noindent\textbf{ACM CCS Concepts.} Theory of computation: Computational geometry; Mathematics of computing: Discrete optimization; Mathematics of computing: Combinatorial optimization; Computing methodologies: Randomized algorithms; Computing methodologies: Genetic algorithms.
\newpage
\section{Introduction}

For a finite point set $P\subset\mathbb{R}^2$, let
\[
  U(P)=\#\bigl\{\{x,y\}\subset P:\|x-y\|=1\bigr\},
\]
and define
\[
  u(n)=\max_{|P|=n} U(P).
\]
The classical unit-distance problem asks for the asymptotic growth of $u(n)$.  Erd\H{o}s conjectured that lattice-type constructions are essentially optimal, in the sense that $u(n)$ should be bounded by $n^{1+o(1)}$ \cite{Erdos1946,OpenAI2026}.  The best general upper bounds remain much larger, of order $O(n^{4/3})$, through incidence-geometric methods related to the Szemer\'edi-Trotter theorem \cite{ST1983,SST1984}.

In 2026, OpenAI announced a counterexample to the Erd\H{o}s conjecture, and a human-verified exposition by Alon, Bloom, Gowers, Litt, Sawin, Shankar, Tsimerman, Wang, and Matchett Wood described the structure of the argument and its relation to earlier number-theoretic ideas \cite{OpenAI2026,AlonEtAl2026}.  Sawin then gave an explicit quantitative refinement, proving that there are arbitrarily large $n$-point sets with more than $n^{1.014}$ unit-distance pairs \cite{Sawin2026}.

As of version~4 of this manuscript, the author was made aware of an active MathOverflow discussion devoted to explicit improvements of the same unit-distance exponent problem \cite{MathOverflowUnitDistance}.  In the same period, Tseng deposited on Zenodo an independently developed certificate package and verification pipeline, reporting the explicit exponent $1+\delta=1.03158935$ and providing machine-checkable data and code for reproducing the certificate \cite{Tseng2026Zenodo}.  In particular, Naslund's MathOverflow contribution reports $\delta>0.03583\ldots$ by combining improved arithmetic and geometric estimates with a reoptimization of the ramified prime set $T$ \cite{NaslundMO2026}; later June 2026 updates in the same discussion report bounds beyond $\delta>0.036$.  These developments were brought to the author's attention after publication of version~3 of this report on arXiv.  It should be noted that some of the newer certificates may not correspond exactly to the optimization problem studied in the present paper: while extending the prime range $T$ appears to be important, some improvements may also involve modified constraint systems or additional degrees of freedom that deviate from Sawin's original formulation.  The results in the present paper are based solely on Sawin's 2026 paper and on the optimization problem explicitly stated there, together with direct extensions of the prime range $T$; the MathOverflow and Zenodo developments are therefore cited as recent related work rather than incorporated into the verified certificates below.

The purpose of this note is narrower than reproving the counterexample.  It studies the finite parameter-selection problem exposed by Sawin's explicit criterion, with an emphasis on reproducible computation and internal validation against the published explicit example.  The central contribution is to formulate the optimization of an improved lower-bound certificate as a nonlinear integer optimization problem over finite data $(T,S_Q,k,R)$, with the real parameter $R$ represented rationally in the implementation, and to apply computationally lightweight but modern integer-optimization heuristics to this formulation.  In particular, the report compares a deterministic greedy construction with a Tailored Integer Evolution Strategy: an integer evolution-strategy variant, in the Rechenberg--Schwefel tradition and following Rudolph's integer-programming mutation model, augmented here with repair operators for the number-theoretic certificate constraints \cite{Rudolph1994}.

The computational contribution is accompanied by a verification pipeline.  The pipeline is first validated by reproducing Sawin's published $n^{1.014\ldots}$ certificate and is then applied to the improved candidates found by the optimization procedures.  The best verified certificate in the present computations supports the cautious clean statement $u(n)>n^{1.0152}$ for arbitrarily large $n$, conditional on Sawin's criterion being applied exactly as cited.  Relatedly, Tseng independently developed a Zenodo-hosted certificate package and verification pipeline for a stronger reported certificate, emphasizing that reproducible, machine-checkable verification has become a central component of this rapidly evolving line of work \cite{Tseng2026Zenodo}.  For readability, Appendix~\ref{app:algorithmic-details} includes a line-by-line walkthrough of Algorithm~\ref{alg:certificate-verification}, explaining the symbols $T$, $S_Q$, $k(p)$, $R$, the quadratic field $K_T$, the splitting tests, the admissibility witnesses, and the final high-precision evaluation.

The paper is organized as follows.  Section~\ref{sec:lattice-visualization} is expository and introduces the geometry of unit-distance graphs through regular-lattice examples in $[0,10]^2$; it is meant to recall the kind of lattice constructions Erd\H{o}s had in mind and is not part of the new certificate.  Section~\ref{sec:integer-programming} formulates Sawin's finite certificate-selection problem as a nonlinear integer optimization problem.  Section~\ref{sec:heuristics} describes the greedy baseline and the Tailored Integer Evolution Strategy, while deferring parameter-level details and pseudocode to Appendix~\ref{app:algorithmic-details}.  Section~\ref{sec:three-bounds} compares the verified certificates, and Section~\ref{sec:reproducible-implementation} summarizes the GitHub implementation and verification commands.  Section~\ref{sec:outlook} discusses remaining room for improvement, and Section~\ref{sec:conclusion} concludes.  The detailed certificate checks, including the validation against Sawin's published example, are collected in Appendix~\ref{app:verification}.

No claim is made here that a coordinate realization of the optimized candidate has been generated.  The author is not a specialist in algebraic number theory; the computations in this note use Sawin's explicit preprint as the mathematical basis for a parameter-optimization and verification study.

\section{Visualization and Solutions on Regular Lattices}\label{sec:lattice-visualization}

This section records two explicit coordinate examples in the fixed square $[0,10]^2$. The idea is to strengthen the intuition about the problem and the kind of solutions that were considered before the 2026 breakthrough, but the section is not essential in terms of results. They are included to clarify the geometry of unit-distance graphs and to provide a concrete comparison with classical lattice constructions.  They are not coordinate realizations of the optimized Sawin-style certificates studied later in the paper.  In both lattice examples below, disks have radius $1/2$ and red segments are drawn only for pairs of centers at Euclidean distance exactly $1$.

\subsection{A local example illustrating touching, overlap, and separation}

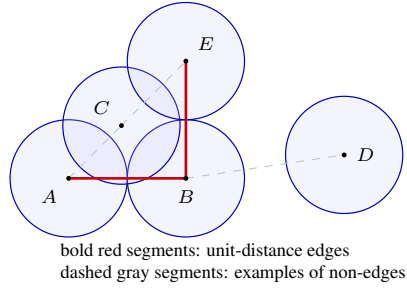
\begin{figure}[htbp]
\centering
\begin{tikzpicture}[
  scale=1.55,
  halfdisk/.style={draw=blue!70!black, fill=blue!20, fill opacity=0.20, line width=0.45pt},
  unitedge/.style={draw=red!80!black, line width=1.0pt},
  pt/.style={circle, fill=black, inner sep=0.65pt},
  faint/.style={draw=gray!50, dashed, line width=0.35pt},
  lab/.style={font=\scriptsize, fill=white, inner sep=1pt}
]
  \coordinate (A) at (0,0);
  \coordinate (B) at (1,0);
  \coordinate (E) at (1,1);
  \coordinate (C) at (0.45,0.45);
  \coordinate (D) at (2.35,0.20);

  \foreach \P in {A,B,C,D,E} {
    \draw[halfdisk] (\P) circle[radius=0.5];
  }

  \draw[unitedge] (A) -- (B);
  \draw[unitedge] (B) -- (E);

  \draw[faint] (A) -- (C);
  \draw[faint] (C) -- (E);
  \draw[faint] (B) -- (D);

  \node[pt,label=below left:{\scriptsize $A$}] at (A) {};
  \node[pt,label=below:{\scriptsize $B$}] at (B) {};
  \node[pt,label=above right:{\scriptsize $E$}] at (E) {};
  \node[pt,label=above left:{\scriptsize $C$}] at (C) {};
  \node[pt,label=right:{\scriptsize $D$}] at (D) {};


  \node[align=left, font=\scriptsize, anchor=west] at (-0.15,-0.72) {bold red segments: unit-distance edges\\dashed gray segments: examples of non-edges};
\end{tikzpicture}
\caption{A local point set showing the three relevant distance regimes simultaneously.  Bold red segments mark all pairs at distance $1$, hence all edges of the unit-distance graph on this five-point set.  The half-unit disks centered at $A$ and $C$ overlap because $\|A-C\|<1$, while the half-unit disks centered at $B$ and $D$ are separated because $\|B-D\|>1$.  Thus overlap between some half-unit disks is allowed; what matters combinatorially is only which pairs are exactly at distance $1$.}
\label{fig:touch-overlap-separate}
\end{figure}

Before discussing the lattice examples, let us briefly remark that the unit-distance problem is not a packing problem, and it is useful to separate three different geometric relations.  If disks of radius $1/2$ are centered at the points of a planar set, then two disks are tangent precisely when their centers are at distance $1$.  However, the problem does not require all other disks to be disjoint.  In particular, one may simultaneously have point pairs at distance $1$ (the edges of the unit-distance graph), point pairs at distance strictly less than $1$ (hence overlapping half-unit disks), and point pairs at distance greater than $1$ (hence separated half-unit disks).  Figure~\ref{fig:touch-overlap-separate} records these three possibilities in a single local example.

\subsection{Hexagonal packing of half-unit disks}

\begin{figure}[p]
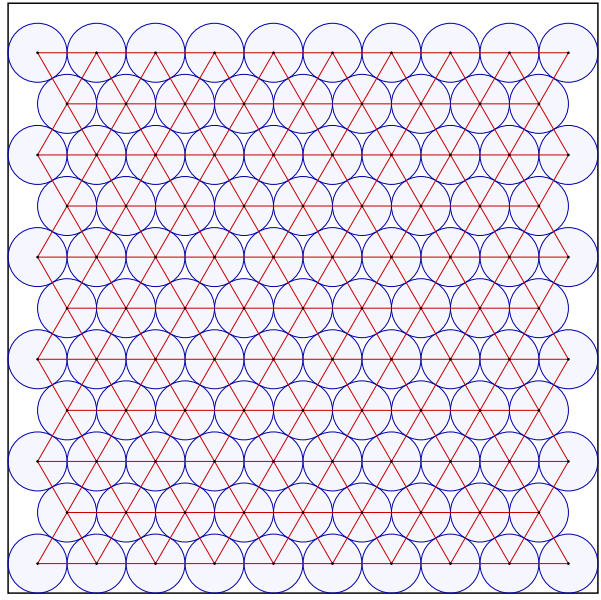

\centering

\caption{Hexagonal packing of half-unit disks in $[0,10]^2$.  There are $105$ centers and $274$ unit-distance pairs.  The disks have radius $1/2$ and are non-overlapping; red segments mark exactly the center pairs at distance $1$.}
\label{fig:hex-half-unit-packing}
\end{figure}

The first example, shown in Figure~\ref{fig:hex-half-unit-packing}, is the standard hexagonal packing of disks of radius $1/2$, clipped so that all disks lie inside $[0,10]^2$; this is the classical optimal planar disk-packing pattern, going back to Fejes T\'oth \cite{FejesToth1942}.  The centers are
\[
  H=\left\{\left(\frac12+i+\frac{\varepsilon_j}{2},\frac12+\frac{\sqrt3}{2}j\right): j=0,\ldots,10\right\},
\]
where $\varepsilon_j=0$ for even rows and $\varepsilon_j=1$ for odd rows.  For even rows one takes $i=0,\ldots,9$, and for odd rows one takes $i=0,\ldots,8$.  Thus there are six rows with ten centers and five rows with nine centers, hence
\[
  |H|=6\cdot 10+5\cdot 9=105.
\]
Unit-distance edges occur horizontally and between adjacent shifted rows.  The horizontal count is
\[
  6\cdot 9+5\cdot 8=94,
\]
and each of the ten adjacent row pairs contributes $18$ diagonal contacts.  Hence
\[
  U(H)=94+10\cdot 18=274.
\]
Equivalently, the average degree of the resulting unit-distance graph is $2U(H)/|H|=548/105\approx5.219$.

\subsection{A multidirectional integer lattice in the same square}

\begin{figure}[p]
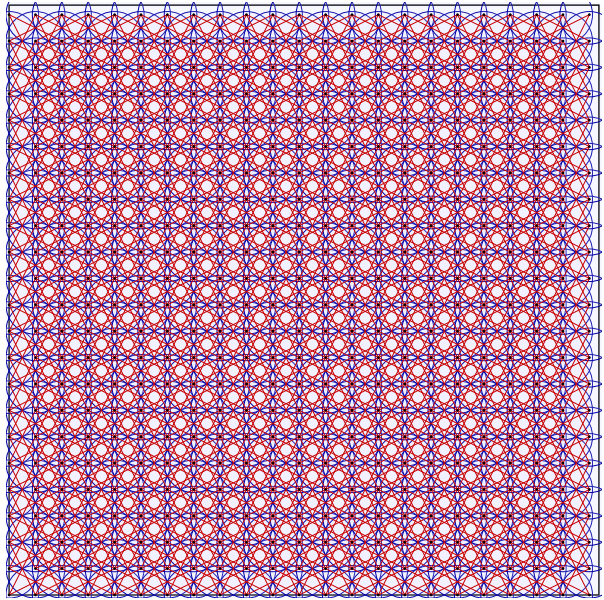

\centering

\caption{Multidirectional lattice arrangement in $[0,10]^2$ based on squared length $5$.  The set has $529$ points and $1848$ unit-distance pairs.  Half-unit disks are shown only as distance markers and may overlap; red segments mark exactly the pairs at distance $1$.}
\label{fig:multidirectional-lattice-0-10}
\end{figure}

The second example uses the integer grid before rescaling, but takes unit distances from several equal-length integer displacement vectors.  Let
\[
  L=\left\{\left(\frac{i}{\sqrt5},\frac{j}{\sqrt5}\right): i,j=0,1,\ldots,22\right\}\subset[0,10]^2.
\]
The upper index $22$ is chosen because $22/\sqrt5<10$, while $23/\sqrt5>10$.  Thus
\[
  |L|=23^2=529.
\]
The four undirected displacement types
\[
  (1,2),\quad (1,-2),\quad (2,1),\quad (2,-1)
\]
have squared length $5$ in the integer grid and therefore become unit displacement vectors after division by $\sqrt5$.  For each type, the number of valid translated copies inside the $23\times23$ integer grid is
\[
  (23-1)(23-2)=22\cdot21.
\]
Consequently
\[
  U(L)=4\cdot22\cdot21=1848.
\]
The average degree is $2U(L)/|L|=3696/529\approx6.988$, which is larger than for the hexagonal packing above.  This comparison is possible because the unit-distance problem is not a packing problem: in Figure~\ref{fig:multidirectional-lattice-0-10}, the half-unit disks are visual distance markers and may overlap, but edges are drawn only between centers exactly one unit apart.

\begin{figure}[htbp]
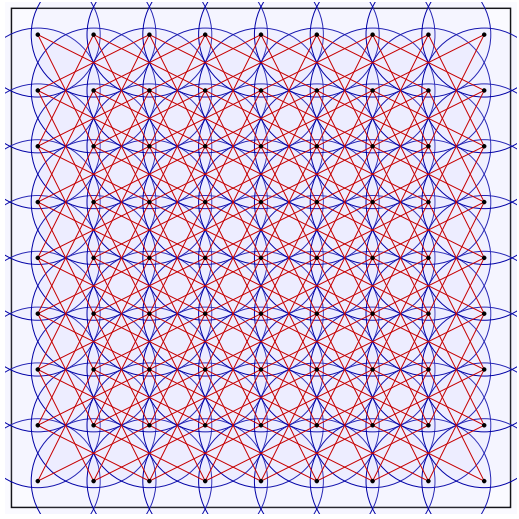

\centering
%
\caption{Local view of the same multidirectional lattice as in Figure~\ref{fig:multidirectional-lattice-0-10}, restricted to the window $[-2,2]^2$.  The points lie on the scaled integer grid $(i/\sqrt5,j/\sqrt5)$, the disks have radius $1/2$, and the red segments show all unit-distance pairs coming from the displacement vectors $(\pm1,\pm2)$ and $(\pm2,\pm1)$.}
\label{fig:multidirectional-lattice-local}
\end{figure}

To make the local structure easier to read, Figure~\ref{fig:multidirectional-lattice-local} also shows a small patch of the same multidirectional lattice in the window $[-2,2]^2$.  This local view makes the repeated pattern generated by the four displacement types $(\pm1,\pm2)$ and $(\pm2,\pm1)$ more apparent, while using exactly the same point set and edge rule as the larger plot.

\section{A nonlinear integer-programming formulation}\label{sec:integer-programming}

The explicit lower-bound criterion suggests the following finite parameter-selection problem \cite{Sawin2026}.  In computational form it is naturally viewed as a nonlinear integer programming problem with arithmetic feasibility constraints.  One chooses:
\begin{itemize}
  \item a finite set $T$ of odd primes;
  \item a finite set $S_Q$ of rational primes;
  \item integer multiplicities $k(p)\geq 1$ for $p\in S_Q$;
  \item a real parameter $R>1$.
\end{itemize}

For the implementation discussed below, define
\[
  e(p)=
  \begin{cases}
  2, & p=2\text{ or }p\in T,\\
  1, & \text{otherwise.}
  \end{cases}
\]
The exponent-gain objective used in the verification is
\begin{align}
\delta(T,S_Q,k,R)&=\frac{N(T,S_Q,k,R)}{D(T,S_Q,k,R)},\label{eq:delta-objective}\\
N(T,S_Q,k,R)&=\log(1-1/R)+\frac12\log(2\pi/e)\notag\\
&\quad+\sum_{p\in S_Q}\frac{1}{4e(p)}\log(k(p)+1)
-\frac18\log\!\left(4\prod_{q\in T}q\right)\notag\\
&\quad-\frac12\log\!\log\!\sqrt{4\prod_{q\in T}q},\notag\\
D(T,S_Q,k,R)&=\log\!\left(2R\prod_{p\in S_Q}p^{k(p)/(2e(p))}+1\right).\notag
\end{align}
The resulting optimization problem is to maximize $\delta(T,S_Q,k,R)$ subject to arithmetic side conditions.  The optimization stage uses a lightweight admissibility model to generate candidates; the candidate reported in Appendix~\ref{app:verification} is then checked separately by exact integer arithmetic for the side conditions implemented in the verification pipeline.

\section{Heuristic optimization strategies}\label{sec:heuristics}

The finite certificate problem can be read as a nonlinear integer optimization problem with arithmetic feasibility constraints.  The variables are the selected prime set $S_Q$, the integer multiplicities $k(p)$, and the rationally represented real parameter $R$.  The role of the optimization algorithms in this paper is to propose candidate certificates; every reported candidate is then checked independently by the verification pipeline described in Appendix~\ref{app:algorithmic-details} and Appendix~\ref{app:verification}.

This section gives only the main algorithmic ideas.  The precise pseudocode, default parameter values, repair rules, and verification loop are collected in Appendix~\ref{app:algorithmic-details}; the numerical certificate checks are recorded in Appendix~\ref{app:verification}.

\subsection{Greedy budget heuristic}\label{subsec:greedy-heuristic}

\begin{figure}[htbp]
\centering
\begin{tikzpicture}[
  node distance=9mm and 9mm,
  box/.style={draw, rounded corners, align=center, minimum height=9mm, minimum width=30mm, fill=blue!6},
  smallbox/.style={draw, rounded corners, align=center, minimum height=8mm, minimum width=29mm, fill=green!8},
  decision/.style={draw, diamond, aspect=2.1, align=center, fill=orange!12, inner sep=1.5pt},
  line/.style={-{Latex[length=2mm]}, thick},
  softline/.style={-{Latex[length=2mm]}, thick, dashed, color=gray!70!black}
]

\node[box] (T) {choose prime set\\$T$};
\node[box, below=of T] (candidates) {enumerate admissible\\candidate primes $p$};
\node[smallbox, below left=of candidates, xshift=-12mm] (benefit) {benefit\\$\frac{1}{4e(p)}\log(k(p)+1)$};
\node[smallbox, below right=of candidates, xshift=12mm] (cost) {cost\\$1$ or $2$};
\node[box, below=20mm of candidates] (score) {rank by\\benefit/cost};
\node[decision, below=of score] (budget) {budget\\available?};
\node[box, below left=of budget, xshift=-11mm] (add) {add $p$\\to $S_Q$};
\node[box, below right=of budget, xshift=11mm] (skip) {skip\\$p$};
\node[box, below=20mm of budget] (scanR) {scan $R$ and evaluate\\$\delta(T,S_Q,k,R)$};
\node[box, below=of scanR] (best) {keep best\\candidate};

\draw[line] (T) -- (candidates);
\draw[line] (candidates) -- (benefit);
\draw[line] (candidates) -- (cost);
\draw[line] (benefit) |- (score);
\draw[line] (cost) |- (score);
\draw[line] (score) -- (budget);
\draw[line] (budget) -- node[left] {yes} (add);
\draw[line] (budget) -- node[right] {no} (skip);
\draw[line] (add.south) |- (scanR.west);
\draw[line] (skip.south) |- (scanR.east);
\draw[line] (scanR) -- (best);

\node[draw, dashed, rounded corners, fit=(benefit)(cost)(score)(budget)(add)(skip), inner sep=5mm] (greedybox) {};
\node[anchor=south west, font=\small, fill=white, inner sep=1pt] at ([xshift=3mm]greedybox.north west) {greedy construction of $S_Q$};
\end{tikzpicture}
\caption{A conceptual optimization pattern for selecting $(T,S_Q,k,R)$.}
\label{fig:optimization-flow}

\end{figure}
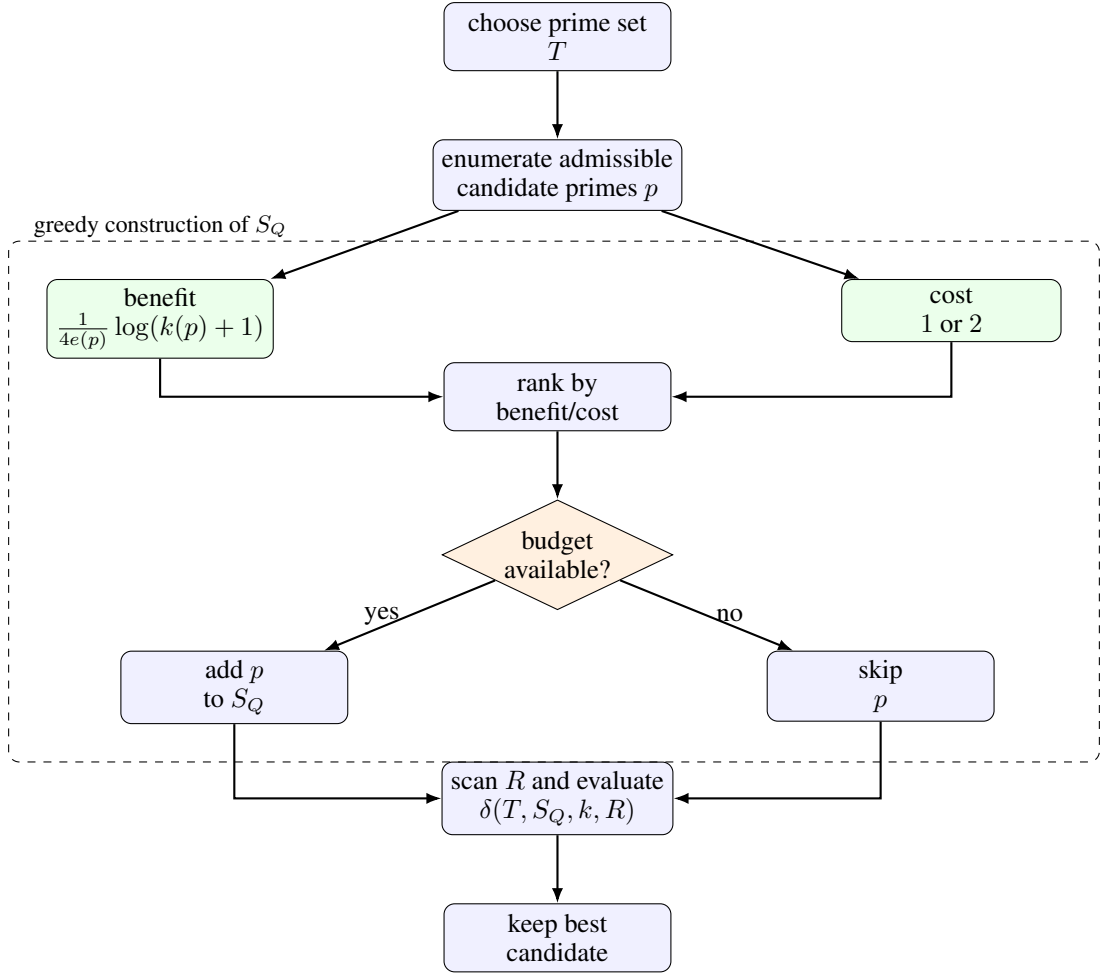

The deterministic greedy baseline starts from the Golod--Shafarevich-type side condition
\[
  \#T+\#S_Q+\#\{p\in S_Q:p\text{ splits in }Q\}+1
  \leq \frac{(\#T-1)^2}{4},
\]
where $Q=\mathbb{Q}\!\left(\sqrt{\prod_{q\in T}q}\right)$ in the simplified quadratic-field picture.  This inequality suggests a knapsack interpretation: a non-splitting prime costs one unit of budget, while a splitting prime costs two.  For a provisional target exponent $c$, the greedy construction assigns the approximate multiplicity
\[
  k(p)\approx \left\lfloor \frac{1}{2c\log p}-1\right\rfloor,
  \qquad
  R\approx 1+\frac{1}{c},
\]
and ranks admissible primes by
\[
  \operatorname{score}(p)=
  \frac{\frac{1}{4e(p)}\log(k(p)+1)}{\operatorname{cost}(p)}.
\]
It then fills the available budget in decreasing score order and finally scans a grid of $R$ values.  Figure~\ref{fig:optimization-flow} summarizes this greedy construction: the dashed frame is the budgeted selection of $S_Q$, while the lower part of the diagram represents the subsequent scan over $R$ and evaluation of $\delta(T,S_Q,k,R)$.

The default greedy run used in the implementation fixes $T$ to Sawin's published prime set, uses candidate primes up to $p_{\max}=300$, sets the provisional target to $c=0.015$, and scans $R\in[40,90]$ on a uniform grid with 500 intervals.  These parameter choices are not part of the theorem; they define a reproducible baseline optimization run.  More detailed algorithmic specifications and verification steps are deferred to Appendix~\ref{app:algorithmic-details}.

\subsection{Tailored Integer Evolution Strategy}\label{subsec:integer-es}

The greedy construction changes the certificate in a largely one-pass way.  The second optimization layer therefore treats the certificate variables directly as integer variables and uses a tailored version of the Integer Evolution Strategy \cite{Rudolph1994} for approximating the solution of the nonlinear integer optimization problem.  This is an integer-valued evolution-strategy instantiation in the Rechenberg--Schwefel tradition.  It uses Rudolph's evolutionary algorithm for integer programming as the relevant integer-mutation model based on the $\ell_1$-symmetric double geometric distribution \cite{Rudolph1994}, and augments it with problem-specific repair operators for the certificate constraints.  Candidate solutions encode selected prime indices, multiplicities, and the numerator $r$ of a rational representation $R=r/s$ with fixed denominator $s$.

The mutation operator is integer-native: it adds two-sided geometric, or discrete-Laplace, steps to integer components instead of mutating real variables and rounding them.  The best reported variant also uses two-parent discrete recombination, where each integer component of an offspring is inherited from one of two parents before mutation.  Since arbitrary mutations can violate the number-theoretic certificate constraints, the implementation augments the integer ES with problem-specific repair operators.  These repairs enforce distinct selected primes, admissible multiplicities, the permitted range for $R$, and the Golod--Shafarevich budget before the objective is evaluated.

The main-text role of this Tailored Integer Evolution Strategy is simple: it optimizes coordinated integer changes that the greedy budget rule is unlikely to find.  The reported single run used a deliberately modest compute budget: population size $\mu=24$, offspring number $\lambda=144$, and $G=160$ generations, corresponding to $24+160\cdot144=23064$ objective-function evaluations, before final verification.  This run completed in less than ten minutes on standard laptop/desktop hardware in the intended reproducibility setting.  Detailed parameters, random seeds, repair rules, and the exact verification loop are specified in Algorithm~\ref{alg:rudolph-discrete-recombination} in Appendix~\ref{app:algorithmic-details}.  The verified certificate produced by this optimization procedure is reported together with the Sawin and greedy certificates in Section~\ref{sec:three-bounds}.

\section{Comparison of verified certificate levels}\label{sec:three-bounds}

The verification pipeline is used in four roles.  It first reproduces Sawin's published explicit certificate as a validation target.  It then checks the greedy certificate, the Tailored integer-ES certificate, and the discrete-recombination variant.  The purpose of this section is only to compare the verified certificate values; the algorithmic details are referenced back to Section~\ref{sec:heuristics} and Appendix~\ref{app:algorithmic-details}.

\begin{table}[htbp]
\centering
\tiny
\begin{tabular}{p{0.28\linewidth}p{0.24\linewidth}p{0.14\linewidth}p{0.20\linewidth}}
\toprule
entry & role & $R$ & exponent parameter $\delta$ \\
\midrule
Sawin published example & validation baseline & $72$ & $0.0141144286784982\ldots$ \\
greedy optimized certificate & deterministic improvement & $66.72240803$ & $0.0151718056372133\ldots$ \\
Tailored integer ES certificate & integer evolutionary improvement & $6672416/100000$ & $0.0152616610684193\ldots$ \\
Tailored integer ES with discrete recombination & recombination variant & $6672416/100000$ & $0.0152628688170072\ldots$ \\
\bottomrule
\end{tabular}
\caption{Verified asymptotic certificate levels. Each row passes the implemented arithmetic checks.}
\label{tab:three-bounds}
\end{table}

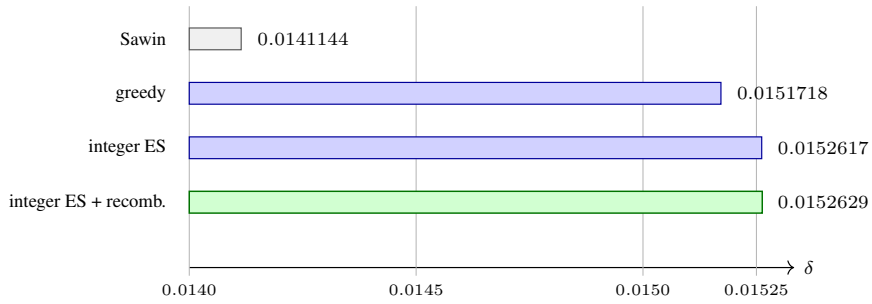
\begin{figure}[htbp]
\centering
\begin{tikzpicture}[
  xscale=1.0, yscale=0.72,
  bar/.style={draw=blue!60!black, fill=blue!18, line width=0.45pt},
  bestbar/.style={draw=green!45!black, fill=green!18, line width=0.55pt},
  base/.style={draw=gray!70!black, fill=gray!12, line width=0.45pt},
  lab/.style={font=\scriptsize, anchor=east},
  val/.style={font=\scriptsize, anchor=west},
  tick/.style={draw=gray!60, line width=0.3pt}
]
  \draw[->, line width=0.45pt] (0,-0.95) -- (8.0,-0.95) node[anchor=west, font=\scriptsize] {$\delta$};
  \foreach \x/\label in {0/0.0140,3/0.0145,6/0.0150,7.5/0.01525} {
    \draw[tick] (\x,-1.10) -- (\x,3.85);
    \node[font=\tiny, anchor=north] at (\x,-1.10) {$\label$};
  }
  \node[lab] at (-0.18,3.25) {Sawin};
  \draw[base] (0,3.05) rectangle (0.6866,3.45);
  \node[val] at (0.78,3.25) {$0.0141144$};

  \node[lab] at (-0.18,2.25) {greedy};
  \draw[bar] (0,2.05) rectangle (7.0308,2.45);
  \node[val] at (7.12,2.25) {$0.0151718$};

  \node[lab] at (-0.18,1.25) {integer ES};
  \draw[bar] (0,1.05) rectangle (7.56997,1.45);
  \node[val] at (7.66,1.25) {$0.0152617$};

  \node[lab] at (-0.18,0.25) {integer ES + recomb.};
  \draw[bestbar] (0,0.05) rectangle (7.57721,0.45);
  \node[val] at (7.66,0.25) {$0.0152629$};

  \node[font=\tiny, anchor=west] at (0,-1.55) {};
\end{tikzpicture}
\caption{Visual comparison of the verified certificate levels.  The horizontal scale is clipped at $\delta=0.014$ so that the differences between the improved certificates remain visible.  The final row is the best certificate found in the recorded computations.}
\label{fig:certificate-comparison}
\end{figure}

The numerical ordering in Table~\ref{tab:three-bounds} and Figure~\ref{fig:certificate-comparison} is the main computational outcome of this version of the report.  Parameter-dashboard visualizations for the Sawin, greedy, and best discrete-recombination certificates are collected in Appendix~\ref{app:verification}; they are displayed in Figures~\ref{fig:sawin-dashboard}, \ref{fig:verified-candidate-dashboard}, and \ref{fig:recomb-candidate-dashboard}.  The greedy optimization already improves the displayed exponent of the validation baseline.  The Tailored Integer Evolution Strategy then improves the greedy certificate by allowing coordinated integer changes in the selected prime set, in the multiplicities $k(p)$, and in the rational representation of $R$.  The discrete-recombination variant gives a further small improvement in this run.  The best current clean consequence supported by the implemented checks is therefore the conservative statement for arbitrarily large $n$:
\[
  u(n)>n^{1.0152}
\]

\section{Reproducible implementation}\label{sec:reproducible-implementation}

The implementations used in this paper are provided in the flat GitHub repository
\begin{center}
\url{https://github.com/emmerichmtm/UnitDistanceProblemOptimizationOfSawinsLowerBound}.
\end{center}
The repository is intended to make the objective function, the verification pipeline, and the optimization heuristics easy to inspect and rerun.  The accompanying README is the authoritative place for file names and command-line usage.  The central verification command is
\begin{verbatim}
python verify_all_certificates_integer_evolution_strategy.py
\end{verbatim}
This script verifies all four certificate levels discussed in the paper: Sawin's published example, the greedy optimized certificate, the Tailored Integer Evolution Strategy certificate, and the discrete-recombination certificate.  It recomputes the finite arithmetic checks for the reported certificates: primality, the parity condition on $T$, Legendre-symbol witnesses, non-splitting of the primes in $S_Q$, the Golod--Shafarevich budget, positivity of the multiplicities, $R>1$, and the value of formula~\eqref{eq:delta-objective}.  The reported output files are
\begin{verbatim}
certificate_verification_all_v23.json
certificate_verification_all_v23.txt
\end{verbatim}
The shorter baseline verifier
\begin{verbatim}
python verify_certificates.py
\end{verbatim}
checks Sawin's published certificate and the greedy optimized certificate and writes
\begin{verbatim}
certificate_verification_results.json
certificate_verification_results.txt
\end{verbatim}

The optimization scripts are deliberately separated.  The deterministic baseline is the greedy strategy of Section~\ref{subsec:greedy-heuristic}:
\begin{verbatim}
python optimize_certificates.py --pmax 300 --c 0.015 --r-min 40 --r-max 90 --r-steps 500
\end{verbatim}
It writes \texttt{optimization\_results.json}.  The improved optimization method is the Tailored Integer Evolution Strategy with discrete recombination from Section~\ref{subsec:integer-es}:
\begin{verbatim}
python rudolph_integer_es_discrete_recombination.py
\end{verbatim}
This script uses \texttt{rudolph\_integer\_ea.py} as a support module and writes
\begin{verbatim}
rudolph_integer_es_discrete_recombination_results.json
rudolph_integer_es_discrete_recombination_results.txt
\end{verbatim}
The reported evolutionary run uses $23064$ objective-function evaluations in a single run and was kept below ten minutes on standard hardware to make reproduction realistic.  Section~\ref{sec:reproducible-implementation} records only the essential commands; full usage notes are kept in the repository README, while Appendix~\ref{app:algorithmic-details} records the algorithmic parameters and Appendix~\ref{app:verification} records the certificate-level results.

For comparison with contemporaneous related work, Tseng's Zenodo record provides an independently developed verification pipeline and certificate package for a stronger reported certificate \cite{Tseng2026Zenodo}.  That implementation was developed independently of the present repository.  Its existence is useful for the field because it offers a separate reproducibility route and a natural target for future cross-verification of certificate formats, arithmetic side conditions, and numerical evaluations.

\section{Outlook and future work}\label{sec:outlook}

The computations reported here suggest several directions for further improvement.  First, the present optimization fixes the same set $T$ used in Sawin's explicit example and explores the induced nonlinear integer-programming problem over $S_Q$, the multiplicities $k(p)$, and a rational representation of $R$.  A natural next step is to enlarge the optimization over $T$ itself, including local exchanges of primes and larger candidate pools.  A first post-v1 improvement in this direction was made available on Zenodo on June 6, 2026: using the optimization and verification pipeline proposed here with an extended prime range, the deposited certificate reports $\delta>0.031$ with $\#T=67$~\cite{EmmerichFrancesco2026Zenodo}.  Second, the greedy and evolutionary optimization pipeline could be complemented by branch-and-bound mixed-integer nonlinear programming heuristics, which have the potential to produce exact optima or upper bounds for quantifying the optimality gap. However, due to the non-linear nature of the objectives and constraints, standard off-the-shelf solvers cannot be directly applied.  Third, the numerical certificate checks should be reproduced with formal interval arithmetic and, e.g., with independent SageMath or Magma verification of the arithmetic side conditions.

There is also a broader theoretical target.  Sawin's paper does not merely provide the explicit lower bound used here; it also discusses limitations and upper bounds for the exponent values obtainable by the same general argument.  The best certificate found in the present computations, $\delta=0.0152628688\ldots$, improves Sawin's displayed explicit value but remains below the upper range discussed for this method.  Thus the computational optimization space has not been exhausted.  Further improvements may come either from better parameter optimization within Sawin's criterion or from strengthening the underlying number-theoretic estimates and field-construction ingredients.

A further direction, brought to the author's attention only after publication of version~3 on arXiv, is to consolidate the optimization and verification pipeline developed here with the newer MathOverflow certificates and the independently developed Zenodo certificate package of Tseng \cite{MathOverflowUnitDistance,Tseng2026Zenodo}.  Tseng's package reports $1+\delta=1.03158935$ together with verification code, while Naslund's contribution and subsequent June 2026 MathOverflow updates indicate that stronger estimates can raise the explicit exponent still further, with reported values above $\delta>0.035$ and beyond $\delta>0.036$ \cite{NaslundMO2026}.  Extending the prime range underlying $T$ appears to be an important ingredient in these improvements.  However, some of the reported certificates may also rely on alternative constraint systems, relaxed conditions, or additional free parameters that deviate from Sawin's original optimization formulation.  The pipeline developed here may therefore be useful not only for reproducing the certificates in this paper, but also as a starting point for independently checking, refining, cross-verifying, and documenting newer certificates in a reproducible form, provided that the verification layer is generalized to cover the corresponding extended formulations.

Finally, the present work remains certificate-level.  Producing a concrete coordinate realization of one of the improved asymptotic certificates would require an explicit finite level of the relevant class-field tower, computable algebraic bases, embeddings, ideals, and a projection/windowing procedure.  That remains a separate and more demanding computational number-theory project.

\section{Conclusion}\label{sec:conclusion}

This report has treated the optimization of improved explicit parameters in Sawin's lower-bound criterion as a nonlinear integer optimization problem.  The optimization variables are the selected primes $S_Q$, their integer multiplicities $k(p)$, and a rational representation of the parameter $R$, subject to arithmetic admissibility and Golod--Shafarevich budget constraints.  The implementation first validates the objective and verification pipeline by reproducing Sawin's published certificate with
\[
  \delta=0.0141144286784982\ldots,
\]
and then applies deterministic greedy optimization and the Tailored Integer Evolution Strategy to the same finite certificate problem.

The best certificate found in the reported computations is obtained by a Tailored Integer Evolution Strategy with discrete recombination.  It has
\[
  \delta=0.0152628688170072\ldots,
\]
passes the implemented arithmetic checks, and supports the clean conservative statement
\[
  u(n)>n^{1.0152}
\]
for arbitrarily large $n$, under the same certificate interpretation used throughout this paper. For extended ranges of $T$, the result obtained with the same optimization-verification pipeline by Emmerich and Cordella after v1 of this article \cite{EmmerichFrancesco2026Zenodo} suggests an even higher value $u(n)>n^{1.031}$ for $\#T = 67$.   These results improve the displayed clean exponent $1.014$ associated with Sawin's explicit parameter choice at the level of verified finite certificate selection.  The result is not a coordinate construction of the corresponding asymptotic planar point sets; the distinction between certificate verification, coordinate realization, and regular-lattice visualizations is maintained throughout the paper.

The main contribution is, therefore, computational and methodological: a state-of-the-art but computationally lightweight integer heuristic, adapted to the nonlinear integer structure of Sawin's certificate model, can find a better explicit certificate than the initially hand-selected one.  The accompanying verification scripts make the certificate checks reproducible, while the regular-lattice figures serve only as explanatory visual material for the unit-distance relation and should not be confused with the algebraic certificate construction.

There remains clear room for further optimization.  More concretely, Proposition~15 in \cite{Sawin2026} gives an abstract upper bound of
\[
  1+\frac{1}{4.116}=1.24295\ldots
\]
for the exponent obtainable from the formulation considered there.  The present certificate is therefore still far below the theoretical ceiling suggested by Sawin's abstraction.

This paper not only establishes a new certificate, but also shows that well-crafted heuristic optimization procedures can contribute to problems in pure mathematics and combinatorial geometry.  Having said so, the paper should also be seen as an invitation to optimize better bounds for this problem, and the openly available implementations of the objective function and verification pipeline can be readily used for this.

\appendix
\section{Pseudocode for the Tailored Integer Evolution Strategy and verification pipeline}\label{app:algorithmic-details}

This appendix gives a precise algorithmic description of the Tailored Integer Evolution Strategy used in the implementation.  The notation is adapted to the certificate problem studied in the report.  A chromosome encodes a finite set of selected prime indices, a vector of integer multiplicities, and an integer numerator for the rational representation of $R$.  The denominator of $R$ is fixed, so every optimization variable is integer-valued.  The objective value is the exponent gain $\delta$ from formula~\eqref{eq:delta-objective}, with infeasible certificates receiving value $-\infty$ or being repaired before evaluation.  The values listed in Algorithm~\ref{alg:rudolph-discrete-recombination} are the default parameters used for the reported run.  The initial population contains the published Sawin certificate, the greedy certificate, the best simple integer-ES certificate, and mutated greedy seeds until the population size is reached.

\begin{algorithm}[H]
\caption{Tailored Integer Evolution Strategy with discrete recombination}
\label{alg:rudolph-discrete-recombination}
\begin{algorithmic}[1]
\Require fixed prime set $T$; candidate-prime bound $p_{\max}=300$; target size $m=22$; budget bound $B$; population size $\mu=24$; offspring number $\lambda=144$; generation limit $G=160$; random seeds $12345,54321,98765,20260601$; rational denominator $s=100000$ for $R=r/s$; $R$-range $40\le R\le90$; multiplicity bounds $1\le k(p)\le80$; mutation parameters $\alpha_I=0.42$, $\alpha_k=0.35$, $\alpha_R=0.35$, $p_I=0.18$, $p_k=0.35$, $p_R=0.85$
\Ensure best feasible certificate $(T,S_Q,k,R)$ and verified value $\delta$
\State Initialize the pseudorandom generator with the selected seed and initialize a population $P$ of $\mu$ feasible integer chromosomes.
\State Each chromosome has the form $x=(I,k,r)$, where $I$ is a list of $m$ indices into $C$, $k$ is the integer multiplicity vector, and $R=r/s$.
\State Evaluate every $x\in P$ by calling the certificate verifier; store $f(x)=\delta(x)$ if feasible and $f(x)=-\infty$ otherwise.
\For{$g=1,2,\ldots,G$}
  \State Let $E$ be the current elite set, consisting of the best feasible parents in $P$.
  \State Initialize an empty offspring set $O$; use the integer mutation parameters $(\alpha_I,\alpha_k,\alpha_R,p_I,p_k,p_R)$ throughout this generation.
  \For{$j=1,2,\ldots,\lambda$}
    \State Select two parents $x^{(a)}=(I^{(a)},k^{(a)},r^{(a)})$ and $x^{(b)}=(I^{(b)},k^{(b)},r^{(b)})$ from $E$.
    \State \textbf{Discrete recombination:} for each integer component, independently inherit the component from parent $a$
    \State $\ $ or parent $b$ with probability $1/2$.
    \State Let the recombined chromosome be $y=(I^{(y)},k^{(y)},r^{(y)})$.
    \State \textbf{Integer-native mutation:} with probability $p_{\rm mut}$ per component, add a two-sided geometric, equivalently \State $\ $ discrete-Laplace, integer step $Z$ to that component.
    \State Repair $y$ by enforcing valid prime indices, distinct selected primes, positive multiplicities, $R>1$, and the \State $\ $Golod--Shafarevich budget.
    \State Evaluate the repaired offspring $y$ using the certificate verifier.
    \State Add $y$ to $O$.
  \EndFor
  \State Set $P$ to the best $\mu$ chromosomes from $P\cup O$ under the objective $f$.
  \State Record the best value at sparse epochs for reproducibility.
\EndFor
\State \Return the best feasible chromosome in the final population.
\end{algorithmic}
\end{algorithm}

\paragraph{Walkthrough.}
The reported run uses $p_{\max}=300$, $m=22$, $\mu=24$, $\lambda=144$, $G=160$, $s=100000$, $40\le R\le90$, and the four random seeds listed above.  Thus a single evolutionary run uses $\mu+G\lambda=24+160\cdot144=23064$ objective-function evaluations, followed by the final verification report.  The optimization starts from feasible certificates, including the published Sawin certificate, the greedy certificate, the best simple integer-ES certificate, and small perturbations of these known good points.  Discrete recombination is used only on integer-coded components: selected-prime indices, multiplicities, and the numerator $r$ in $R=r/s$.  After recombination, integer-valued mutation in the Rudolph tradition changes integer components directly; it does not mutate real-valued variables and then round them.  A repair step keeps the candidate within the admissible optimization space.  Finally, every survivor is checked by the same verification routine used for the reported certificates.  This makes the optimization layer heuristic, but the scoring layer deterministic and reproducible.

The mutation law is the essential Rudolph-style ingredient: the step $Z\in\mathbb Z$ is $\ell_1$ symmetric, centered at zero, and has geometrically decaying tails.  Thus small moves are most frequent, while larger jumps remain possible.  This is appropriate for the present nonlinear integer-programming variant, where changing an optimization variable by a small amount can alter the certificate value, but occasional larger moves help escape local plateaus.

\begin{algorithm}[H]
\caption{Certificate verification pipeline}
\label{alg:certificate-verification}
\begin{algorithmic}[1]
\Require candidate data $(T,S_Q,k,R)$
\Ensure pass/fail status and exponent value $\delta$ if all checks pass
\State Check that $T$ consists of odd primes and that the required parity condition on primes $q\equiv 3\pmod 4$ is satisfied.
\State Compute $A_T=\prod_{q\in T}q$ and set $K_T=\mathbb Q(\sqrt{A_T})$ at the level of exact integer arithmetic.
\For{each $p\in S_Q$}
  \State Check primality of $p$ and positivity of $k(p)$.
  \State Compute the Legendre/Kronecker symbols needed to determine whether $p$ is ramified, inert, or split in $K_T$.
  \State Check the admissibility condition: $p\equiv 1\pmod 4$ or $p$ is inert in $\mathbb Q(\sqrt q)$ for at least one $q\in T$.
\EndFor
\State Compute the split-prime count $s_Q=\#\{p\in S_Q:p\text{ splits in }K_T\}$.
\State Verify the budget inequality
\[
\#T+\#S_Q+s_Q+1\leq \frac{(\#T-1)^2}{4}.
\]
\State Check $R>1$ and evaluate $N(T,S_Q,k,R)$ and $D(T,S_Q,k,R)$ in formula~\eqref{eq:delta-objective} using high-precision decimal arithmetic.
\State Return $\delta=N/D$ together with all diagnostic flags.
\end{algorithmic}
\end{algorithm}

\paragraph{Detailed walkthrough of Algorithm~\ref{alg:certificate-verification}.}
The verifier is a computational certificate checker based on the theory outlined in Sawin \cite{Sawin2026} that can be independent of the optimizer used.  Its input is a finite certificate
\begin{equation*}
  (T,S_Q,k,R),
\end{equation*}
where $T$ is the prime set defining the quadratic field, $S_Q$ is the selected auxiliary prime set, $k:S_Q\to\mathbb Z_{>0}$ assigns a positive integer multiplicity to every selected prime, and $R$ is the positive real parameter appearing in the exponent formula.  The subscript in $S_Q$ is only a name for the selected set attached to the quadratic-field condition; it should not be confused with the rational field $\mathbb Q$.

The first line checks the shape of $T$.  Every element of $T$ must be an ordinary rational prime, hence a prime number in $\mathbb Z$, and all of these primes must be odd.  In addition, the number of primes $q\in T$ with $q\equiv 3\pmod 4$ must be odd.  This parity condition fixes the congruence class of the squarefree radicand used below and keeps the implemented certificate in the case for which the stated formula and side conditions are being applied.

The second line constructs the quadratic field used for the splitting tests.  The verifier forms
\begin{equation*}
  A_T=\prod_{q\in T}q,
  \qquad
  K_T=\mathbb Q(\sqrt{A_T}).
\end{equation*}
The notation $A_T$ is used here for the radicand in order not to confuse it with the denominator $D(T,S_Q,k,R)$ in the objective formula~\eqref{eq:delta-objective}.  Since $T$ is a set of distinct primes, $A_T$ is squarefree.  All arithmetic used to build $A_T$ and to test divisibility and congruences is exact integer arithmetic.

The loop over $p\in S_Q$ then checks the selected primes one by one.  First, $p$ must itself be a rational prime.  Second, its multiplicity $k(p)$ must be a positive integer.  The multiplicity is the exponent-like weight with which $p$ contributes to the denominator product and to the logarithmic gain term in formula~\eqref{eq:delta-objective}.  The auxiliary function
\begin{equation*}
  e(p)=
  \begin{cases}
  2, & p=2\text{ or }p\in T,\\
  1, & \text{otherwise}
  \end{cases}
\end{equation*}
records the local normalization used in that formula.  Thus selected primes that also lie in $T$, and the prime $2$, are treated with $e(p)=2$.

The next check determines whether a selected prime is split, inert, or ramified in $K_T$.  For an odd prime $p\nmid A_T$, this is governed by the Legendre symbol
\begin{equation*}
  \chi_T(p)=\left(\frac{A_T}{p}\right).
\end{equation*}
This symbol is not a fraction.  It is equal to $1$ if $A_T$ is a square modulo $p$, equal to $-1$ if $A_T$ is not a square modulo $p$, and equal to $0$ if $p\mid A_T$.  Consequently,
\begin{align*}
  \chi_T(p)=1  &\quad\Longleftrightarrow\quad p\text{ splits in }K_T,\\
  \chi_T(p)=-1 &\quad\Longleftrightarrow\quad p\text{ is inert in }K_T,\\
  \chi_T(p)=0  &\quad\Longleftrightarrow\quad p\text{ ramifies in }K_T.
\end{align*}
For $p=2$, the implementation uses the corresponding Kronecker-symbol convention.  In the reported certificates, selected primes are intended not to split in $K_T$; equivalently, the strict desired outcome is that no selected $p$ has $\chi_T(p)=1$.  The verifier records the split count $s_Q$ explicitly so that the Golod--Shafarevich budget line can be checked transparently.

A useful toy example is $T=\{3,7,11\}$, for which $A_T=231$ and $K_T=\mathbb Q(\sqrt{231})$.  For $p=17$, one has $231\equiv 10\pmod {17}$, and $10$ is not a quadratic residue modulo $17$, so $\left(\frac{231}{17}\right)=-1$ and $17$ is inert.  For $p=5$, one has $231\equiv1\pmod5$, so $\left(\frac{231}{5}\right)=1$ and $5$ splits.  For $p=7$, one has $7\mid231$, so $\left(\frac{231}{7}\right)=0$ and $7$ ramifies.

The admissibility test is separate from the nonsplitting test in $K_T$.  In the implemented form of the criterion, a selected prime $p\in S_Q$ is accepted if either
\begin{equation*}
  p\equiv1\pmod4,
\end{equation*}
or there is at least one prime $q\in T$ for which $p$ is inert in the simpler quadratic field $\mathbb Q(\sqrt q)$.  For odd $p\nmid q$, the latter condition means
\begin{equation*}
  \left(\frac{q}{p}\right)=-1.
\end{equation*}
An admissibility witness is the concrete data proving this step.  If $p\equiv1\pmod4$, the congruence itself is the witness.  Otherwise, a witness is a specific prime $q\in T$ such that $\left(\frac{q}{p}\right)=-1$; for $p=2$ the corresponding Kronecker-symbol test is used.  The witness tables in Appendix~\ref{app:verification} list exactly such certificates, so the check is reproducible prime by prime.

After the loop, the verifier evaluates the Golod--Shafarevich budget
\begin{equation*}
  \#T+\#S_Q+s_Q+1\leq \frac{(\#T-1)^2}{4}.
\end{equation*}
The term $\#T$ counts the primes used to define $K_T$, the term $\#S_Q$ counts the selected auxiliary primes, $s_Q$ is the number of selected primes that split in $K_T$, and the final $+1$ is the constant contribution appearing in this simplified budget.  Thus a nonsplitting selected prime costs one unit of budget, while a splitting selected prime would cost an additional unit.  The verified certificates in this report are arranged so that $s_Q=0$, hence the no-splitting condition and the budget check agree cleanly.

The last line checks $R>1$ and evaluates the objective formula~\eqref{eq:delta-objective}.  In that formula, the numerator $N(T,S_Q,k,R)$ contains the logarithmic gain terms, the discriminant-like penalty from $4A_T$, and the logarithmic correction term, while the denominator $D(T,S_Q,k,R)$ contains the logarithm of the size parameter
\begin{equation*}
  2R\prod_{p\in S_Q}p^{k(p)/(2e(p))}+1.
\end{equation*}
The verifier uses high-precision decimal arithmetic for these transcendental evaluations.  This is sufficient for reproducible numerical checking of the displayed decimals; a fully formal proof certificate would ideally replace the final floating-point step by interval or rationally certified bounds.

The public implementation of this pipeline is available at
\begin{center}
\url{https://github.com/emmerichmtm/UnitDistanceProblemOptimizationOfSawinsLowerBound}
\end{center}
Algorithm~\ref{alg:certificate-verification} summarizes the verification pipeline.  The verification scripts are meant to be readable reference implementations: they reproduce Sawin's published certificate, the greedy certificate, the Tailored Integer Evolution Strategy certificate, and the discrete-recombination certificate reported in this paper.

\section{Validation of the verification pipeline and verification of the best-found candidate}\label{app:verification}

This appendix records four verification examples based on the explicit criterion used in Sawin's quantitative lower bound \cite{Sawin2026}.  First, the pipeline is run on Sawin's published parameter choice from the proof of Theorem~1, in order to verify that the implementation reproduces the advertised exponent $\delta\approx 0.0141144$.  Second, the same pipeline is applied to the greedy certificate produced by the deterministic optimization procedure.  Third, it is applied to the certificate produced by the Tailored Integer Evolution Strategy. Fourth, it is applied to the certificate produced by the discrete-recombination variant.  In all cases, the checked data consist of finite sets $T$ and $S_Q$, integer weights $k(p)$, and a real parameter $R>1$, from which the exponent in formula~\eqref{eq:delta-objective} is computed once the arithmetic side conditions are verified.  The OpenAI announcement and the human-verified remarks explain the broader counterexample and context \cite{OpenAI2026,AlonEtAl2026}.

\subsection*{Validation on Sawin's published parameter choice}

The published example in the proof of Theorem~1 of \cite{Sawin2026} uses
\[
T=\{3,5,7,11,13,17,19,23,29,31,37,41,43\},
\]
with
\[
\begin{split}
S_Q=\{&2,3,5,7,11,13,17,19,23,29,47,71,79,97,101,107,109,139,151,163,167,179\}.
\end{split}
\]
The multiplicities are
\[
\begin{array}{c|rrrrrrrrrrr}
p&2&3&5&7&11&13&17&19&23&29&47\\
\hline
k(p)&50&31&21&17&14&13&12&11&10&10&8
\end{array}
\]
\[
\begin{array}{c|rrrrrrrrrrr}
p&71&79&97&101&107&109&139&151&163&167&179\\
\hline
k(p)&7&7&7&7&7&7&6&6&6&6&6
\end{array}
\]
and
\[
R=72.
\]
The purpose of including this example here is twofold: it validates the implementation against the published reference point, and it illustrates the certificate language before the improved candidate is discussed.

The same arithmetic checks as in the optimized case pass for this published example.  The primes in $T$ are odd, exactly seven of them are congruent to $3\pmod 4$, the set $S_Q$ contains no split primes in $Q=\mathbb{Q}(\sqrt D)$ with $D=\prod_{q\in T}q$, and the Golod--Shafarevich budget is again exactly saturated:
\[
\#T+\#S_Q+\#\{p\in S_Q:p\text{ splits in }Q\}+1=13+22+0+1=36=\frac{(13-1)^2}{4}.
\]
Evaluating formula~\eqref{eq:delta-objective} gives
\[
\text{numerator}=3.8822487482003876\ldots,
\qquad
\text{denominator}=275.0553236430010\ldots,
\]
so that
\[
\delta=0.014114428678498239\ldots,
\]
which reproduces Sawin's printed value $0.014114\ldots$.

\begin{proposition}[Validation against Sawin's published example]
Using the finite data displayed above, the verification pipeline reproduces the explicit lower-bound certificate from the proof of Theorem~1 in \cite{Sawin2026} and confirms the numerical exponent $\delta\approx 0.0141144287$.
\end{proposition}

\subsection*{Graphical description of Sawin's published example}

Figure~\ref{fig:sawin-dashboard} visualizes the published certificate at the parameter level.  It records the set $T$, the selected primes $S_Q$, their multiplicities $k(p)$, and the exact budget saturation.

\begin{figure}[htbp]
\centering
\resizebox{0.98\linewidth}{!}{%
\begin{tikzpicture}[
  font=\scriptsize,
  tcell/.style={draw=blue!70!black, fill=blue!9, rounded corners=1.5pt, minimum width=0.70cm, minimum height=0.47cm, align=center},
  inertcell/.style={draw=green!45!black, fill=green!10, rounded corners=1.5pt, minimum width=0.90cm, minimum height=0.70cm, align=center},
  ramcell/.style={draw=orange!75!black, fill=orange!18, rounded corners=1.5pt, minimum width=0.90cm, minimum height=0.70cm, align=center},
  basecell/.style={draw=gray!70!black, fill=gray!15, rounded corners=1.5pt, minimum width=0.35cm, minimum height=0.25cm, align=center},
  note/.style={align=left, font=\scriptsize},
  lab/.style={font=\scriptsize\bfseries, anchor=east}
]

\node[font=\small\bfseries, align=center] at (5.55,1.05) {Published Sawin-example fingerprint};
\node[note, align=center] at (5.55,0.55) {$Q=\mathbb{Q}(\sqrt{D})$, $D=\prod_{q\in T}q=6541380665835015$, $R=72$, $\delta=0.0141144287\ldots$};

\node[lab] at (-0.25,0) {$T$};
\foreach \p [count=\i from 0] in {3,5,7,11,13,17,19,23,29,31,37,41,43} {
  \node[tcell] at (\i*0.88,0) {\p};
}
\node[note, anchor=west] at (0,-0.55) {13 odd primes; exactly seven are congruent to $3\bmod 4$.};

\node[lab] at (-0.25,-1.35) {$S_Q$};
\foreach \p/\k/\sty [count=\i from 0] in {2/50/ramcell,3/31/ramcell,5/21/ramcell,7/17/ramcell,11/14/ramcell,13/13/ramcell,17/12/ramcell,19/11/ramcell,23/10/ramcell,29/10/ramcell,47/8/inertcell} {
  \node[\sty] at (\i*1.02,-1.35) {\textbf{\p}\\$k=\k$};
}
\foreach \p/\k/\sty [count=\i from 0] in {71/7/inertcell,79/7/inertcell,97/7/inertcell,101/7/inertcell,107/7/inertcell,109/7/inertcell,139/6/inertcell,151/6/inertcell,163/6/inertcell,167/6/inertcell,179/6/inertcell} {
  \node[\sty] at (\i*1.02,-2.32) {\textbf{\p}\\$k=\k$};
}
\node[note, anchor=west] at (0,-2.95) {22 selected primes; the first 10 are ramified in $Q$, the last 12 are inert in $Q$; no selected prime splits in $Q$.};

\node[lab] at (-0.25,-3.85) {budget};
\foreach \i in {0,...,12} {\draw[draw=blue!60!black, fill=blue!20] ({0.285*\i},-4.00) rectangle ++(0.265,0.32);}
\foreach \i in {13,...,34} {\draw[draw=green!45!black, fill=green!18] ({0.285*\i},-4.00) rectangle ++(0.265,0.32);}
\draw[draw=gray!70!black, fill=gray!20] ({0.285*35},-4.00) rectangle ++(0.265,0.32);
\draw[decorate, decoration={brace, amplitude=3pt}, thick] (0,-4.25) -- ({0.285*13-0.02},-4.25) node[midway, yshift=-0.36cm] {$\#T=13$};
\draw[decorate, decoration={brace, amplitude=3pt}, thick] ({0.285*13},-4.25) -- ({0.285*35-0.02},-4.25) node[midway, yshift=-0.36cm] {$\#S_Q=22$};
\draw[decorate, decoration={brace, amplitude=3pt}, thick] ({0.285*35},-4.25) -- ({0.285*36-0.02},-4.25) node[midway, yshift=-0.36cm] {$+1$};
\node[note, anchor=west] at (0,-4.95) {$13+22+0+1=36=((13-1)^2)/4$: the Golod--Shafarevich budget is exactly saturated.};

\node[ramcell, minimum width=0.35cm, minimum height=0.25cm] at (0,-5.65) {};
\node[anchor=west] at (0.27,-5.65) {ramified in $Q$};
\node[inertcell, minimum width=0.35cm, minimum height=0.25cm] at (2.15,-5.65) {};
\node[anchor=west] at (2.42,-5.65) {inert in $Q$};
\node[tcell, minimum width=0.35cm, minimum height=0.25cm] at (4.05,-5.65) {};
\node[anchor=west] at (4.32,-5.65) {prime in $T$};
\node[basecell] at (5.80,-5.65) {};
\node[anchor=west] at (6.07,-5.65) {constant $+1$ budget term};
\end{tikzpicture}}
\caption{Parameter-level visualization of Sawin's published explicit example.  The figure is included here as a validation target for the verification pipeline and as a certificate-level illustration of the method.}
\label{fig:sawin-dashboard}
\end{figure}

\subsection*{Greedy optimized candidate from the deterministic optimization}

After the validation run above, the same pipeline was applied to the best candidate found by the optimization procedure.  That candidate is
\[
T=\{3,5,7,11,13,17,19,23,29,31,37,41,43\},
\]
with
\[
\begin{split}
S_Q=\{&2,3,47,71,79,97,101,107,109,139,151,163,167,179,\\
      &191,211,223,239,241,251,257,263\}.
\end{split}
\]
The chosen multiplicities are
\[
\begin{array}{c|rrrrrrrrrrr}
p&2&3&47&71&79&97&101&107&109&139&151\\
\hline
k(p)&47&29&7&6&6&6&6&6&6&5&5
\end{array}
\]
\[
\begin{array}{c|rrrrrrrrrrr}
p&163&167&179&191&211&223&239&241&251&257&263\\
\hline
k(p)&5&5&5&5&5&5&5&5&5&5&4
\end{array}
\]
and the grid-selected value is
\[
R=66.72240803.
\]

\subsection*{Visualization of the verified certificate}

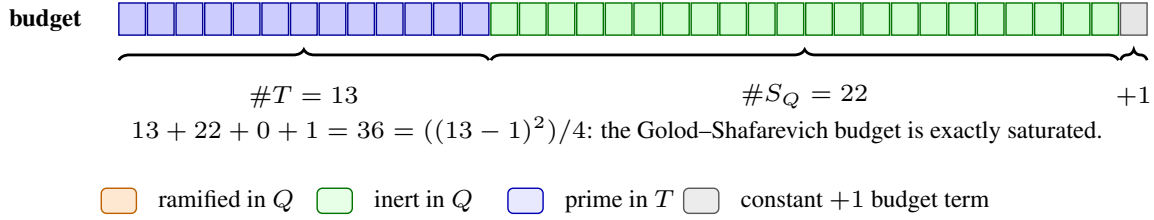
\begin{figure}[htbp]
\centering
\resizebox{0.98\linewidth}{!}{%
\begin{tikzpicture}[
  font=\scriptsize,
  tcell/.style={draw=blue!70!black, fill=blue!9, rounded corners=1.5pt, minimum width=0.70cm, minimum height=0.47cm, align=center},
  inertcell/.style={draw=green!45!black, fill=green!10, rounded corners=1.5pt, minimum width=0.90cm, minimum height=0.70cm, align=center},
  ramcell/.style={draw=orange!75!black, fill=orange!18, rounded corners=1.5pt, minimum width=0.90cm, minimum height=0.70cm, align=center},
  basecell/.style={draw=gray!70!black, fill=gray!15, rounded corners=1.5pt, minimum width=0.35cm, minimum height=0.25cm, align=center},
  note/.style={align=left, font=\scriptsize},
  lab/.style={font=\scriptsize\bfseries, anchor=east}
]

\node[font=\small\bfseries, align=center] at (5.55,1.05) {Optimized candidate fingerprint};
\node[note, align=center] at (5.55,0.55) {$Q=\mathbb{Q}(\sqrt{D})$, $D=\prod_{q\in T}q=6541380665835015$, $R=66.72240803$, $\delta=0.0151718056\ldots$};

\node[lab] at (-0.25,0) {$T$};
\foreach \p [count=\i from 0] in {3,5,7,11,13,17,19,23,29,31,37,41,43} {
  \node[tcell] at (\i*0.88,0) {\p};
}
\node[note, anchor=west] at (0,-0.55) {13 odd primes; exactly seven are congruent to $3\bmod 4$.};

\node[lab] at (-0.25,-1.35) {$S_Q$};
\foreach \p/\k/\sty [count=\i from 0] in {2/47/ramcell,3/29/ramcell,47/7/inertcell,71/6/inertcell,79/6/inertcell,97/6/inertcell,101/6/inertcell,107/6/inertcell,109/6/inertcell,139/5/inertcell,151/5/inertcell} {
  \node[\sty] at (\i*1.02,-1.35) {\textbf{\p}\\$k=\k$};
}
\foreach \p/\k/\sty [count=\i from 0] in {163/5/inertcell,167/5/inertcell,179/5/inertcell,191/5/inertcell,211/5/inertcell,223/5/inertcell,239/5/inertcell,241/5/inertcell,251/5/inertcell,257/5/inertcell,263/4/inertcell} {
  \node[\sty] at (\i*1.02,-2.32) {\textbf{\p}\\$k=\k$};
}
\node[note, anchor=west] at (0,-2.95) {22 selected primes; orange means ramified in $Q$, green means inert in $Q$; no selected prime splits in $Q$.};

\node[lab] at (-0.25,-3.85) {budget};
\foreach \i in {0,...,12} {
  \draw[draw=blue!60!black, fill=blue!20] ({0.285*\i},-4.00) rectangle ++(0.265,0.32);
}
\foreach \i in {13,...,34} {
  \draw[draw=green!45!black, fill=green!18] ({0.285*\i},-4.00) rectangle ++(0.265,0.32);
}
\draw[draw=gray!70!black, fill=gray!20] ({0.285*35},-4.00) rectangle ++(0.265,0.32);
\draw[decorate, decoration={brace, amplitude=3pt}, thick] (0,-4.25) -- ({0.285*13-0.02},-4.25) node[midway, yshift=-0.36cm] {$\#T=13$};
\draw[decorate, decoration={brace, amplitude=3pt}, thick] ({0.285*13},-4.25) -- ({0.285*35-0.02},-4.25) node[midway, yshift=-0.36cm] {$\#S_Q=22$};
\draw[decorate, decoration={brace, amplitude=3pt}, thick] ({0.285*35},-4.25) -- ({0.285*36-0.02},-4.25) node[midway, yshift=-0.36cm] {$+1$};
\node[note, anchor=west] at (0,-4.95) {$13+22+0+1=36=((13-1)^2)/4$: the Golod--Shafarevich budget is exactly saturated.};

\node[ramcell, minimum width=0.35cm, minimum height=0.25cm] at (0,-5.65) {};
\node[anchor=west] at (0.27,-5.65) {ramified in $Q$};
\node[inertcell, minimum width=0.35cm, minimum height=0.25cm] at (2.15,-5.65) {};
\node[anchor=west] at (2.42,-5.65) {inert in $Q$};
\node[tcell, minimum width=0.35cm, minimum height=0.25cm] at (4.05,-5.65) {};
\node[anchor=west] at (4.32,-5.65) {prime in $T$};
\node[basecell] at (5.80,-5.65) {};
\node[anchor=west] at (6.07,-5.65) {constant $+1$ budget term};
\end{tikzpicture}}
\caption{Parameter-level visualization of the optimized candidate in Appendix~\ref{app:verification}. The colored budget bar shows exact saturation of the side condition: 13 slots for $T$, 22 slots for $S_Q$, no split-prime penalty, and one constant slot.}
\label{fig:verified-candidate-dashboard}
\end{figure}

Figure~\ref{fig:verified-candidate-dashboard} visualizes the optimized candidate at the parameter level.  It is not a drawing of the final Euclidean point configuration, whose construction is number-theoretic and asymptotic; rather, it shows the finite certificate data checked by the implementation: the prime set $T$, the selected primes $S_Q$, their multiplicities $k(p)$, and the exact saturation of the Golod--Shafarevich budget.

\subsection*{Coordinate-realization scope}

The dashboard figures in this appendix visualize finite certificate data, not Euclidean coordinates of the asymptotic point sets.  A literal coordinate realization of an optimized certificate would require additional number-theoretic data not specified by $(T,S_Q,k,R)$: an explicit finite level of the relevant class-field tower, algebraic bases or defining polynomials, ideals, embeddings, normalization, and a bounded planar window.  The present implementation therefore focuses on parameter optimization and certificate verification.  Constructing such coordinates remains a separate computational number-theory project.

\subsection*{Arithmetic checks}

Let
\(
D=\prod_{q\in T}q=6541380665835015,
\qquad
Q=\mathbb{Q}(\sqrt D).
\)
The set $T$ consists of odd primes and exactly seven of its elements are congruent to $3\pmod 4$, so the required parity condition is satisfied.  The verification checks, using exact integer arithmetic and Legendre-symbol computations, that every prime in $S_Q$ is either ramified or inert in $Q$; none is split.  Thus
\[
\#\{p\in S_Q:p\text{ splits in }Q\}=0.
\]
The Golod--Shafarevich budget side condition therefore becomes
\[
\#T+\#S_Q+\#\{p\in S_Q:p\text{ splits in }Q\}+1
=13+22+0+1=36, \mbox{while}
\]
\[
\frac{(\#T-1)^2}{4}=\frac{12^2}{4}=36.
\]
Hence the budget is exactly saturated.

For admissibility in Sawin's Lemma~12, each $p\in S_Q$ must be congruent to $1\pmod 4$ or inert in $\mathbb{Q}(\sqrt q)$ for some $q\in T$.  The verification script found the following witnesses.

\begin{center}
\scriptsize
\begin{tabular}{rll}
\toprule
$p$ & status in $Q$ & admissibility witness \\
\midrule
2   & ramified & inert in $\mathbb{Q}(\sqrt5)$ \\
3   & ramified & inert in $\mathbb{Q}(\sqrt5)$ \\
47  & inert & inert in $\mathbb{Q}(\sqrt5)$ \\
71  & inert & inert in $\mathbb{Q}(\sqrt7)$ \\
79  & inert & inert in $\mathbb{Q}(\sqrt3)$ \\
97  & inert & $p\equiv1\pmod4$ \\
101 & inert & $p\equiv1\pmod4$ \\
107 & inert & inert in $\mathbb{Q}(\sqrt5)$ \\
109 & inert & $p\equiv1\pmod4$ \\
139 & inert & inert in $\mathbb{Q}(\sqrt3)$ \\
151 & inert & inert in $\mathbb{Q}(\sqrt3)$ \\
163 & inert & inert in $\mathbb{Q}(\sqrt3)$ \\
167 & inert & inert in $\mathbb{Q}(\sqrt5)$ \\
179 & inert & inert in $\mathbb{Q}(\sqrt7)$ \\
191 & inert & inert in $\mathbb{Q}(\sqrt7)$ \\
211 & inert & inert in $\mathbb{Q}(\sqrt3)$ \\
223 & inert & inert in $\mathbb{Q}(\sqrt3)$ \\
239 & inert & inert in $\mathbb{Q}(\sqrt7)$ \\
241 & inert & $p\equiv1\pmod4$ \\
251 & inert & inert in $\mathbb{Q}(\sqrt{11})$ \\
257 & inert & $p\equiv1\pmod4$ \\
263 & inert & inert in $\mathbb{Q}(\sqrt5)$ \\
\bottomrule
\end{tabular}
\end{center}

\subsection*{Numerical check}

Using 80-digit decimal arithmetic in the standard-library Python module \texttt{decimal}, the implementation evaluates formula~\eqref{eq:delta-objective} for the optimized candidate as follows:
\[
\begin{aligned}
\text{numerator}   &=4.3342249963460248\ldots,\\
\text{denominator} &=285.6762800674879\ldots,\\
\delta             &=0.01517180563721325\ldots .
\end{aligned}
\]
For the conservative target $0.015$, the margin is
\[
\text{numerator}-0.015\cdot \text{denominator}
=0.04908079533370595\ldots>0.
\]
The full 80-digit values are written to \texttt{verification\_attempt.txt}.
This margin is large compared with ordinary floating-point roundoff, although a publication-grade statement should still use formal interval arithmetic.

\subsection*{Tailored Integer Evolution Strategy certificate}

The best certificate found by the Tailored Integer Evolution Strategy keeps the same set $T$ but changes the selected primes, several multiplicities, and $R$.  The real parameter is encoded rationally as
\[
R=\frac{6672416}{100000}=66.72416.
\]
The selected prime set is
\[
\begin{split}
S_Q=\{&2,3,5,47,71,79,97,101,107,109,139,151,163,167,179,\\
      &191,211,223,239,241,251,263\}.
\end{split}
\]
The multiplicities are
\[
\begin{array}{c|rrrrrrrrrrr}
p&2&3&5&47&71&79&97&101&107&109&139\\
\hline
k(p)&49&29&19&7&7&6&6&6&6&6&6
\end{array}
\]
\[
\begin{array}{c|rrrrrrrrrrr}
p&151&163&167&179&191&211&223&239&241&251&263\\
\hline
k(p)&5&5&6&5&5&5&5&5&5&5&5
\end{array}
\]
The verification checks again pass: the budget is exactly saturated,
\[
13+22+0+1=36=\frac{(13-1)^2}{4},
\]
and no selected prime splits in $Q=\mathbb{Q}(\sqrt D)$.  Evaluating formula~\eqref{eq:delta-objective} gives
\[
\begin{aligned}
\text{numerator}   &=4.4218933893294341\ldots,\\
\text{denominator} &=289.73867061427447\ldots,\\
\delta             &=0.01526166106841929\ldots .
\end{aligned}
\]
For the clean target $0.0152$, the numerical margin is positive:
\[
4.4218933893294341\ldots-0.0152\cdot 289.73867061427447\ldots>0.
\]
Thus this certificate supports the clean exponent $1.0152$, subject to the same mathematical caveats and independent-check requirements as the greedy certificate.

\begin{proposition}[Lower-bound consequence]
Assuming Sawin's explicit criterion is applied exactly as in \cite{Sawin2026}, the Tailored Integer Evolution Strategy certificate above supports the clean bound
\[
  u(n)>n^{1.0152}
\]
for arbitrarily large $n$.
\end{proposition}

\subsection*{Tailored Integer Evolution Strategy with discrete recombination}

A two-parent discrete-recombination variant was also tested.  Each offspring inherits the selected-prime indices, multiplicities, and rational numerator of $R$ componentwise from two selected parents, after which the same integer-native mutation operator is applied.  This variant produced a small additional improvement over the simpler Tailored Integer Evolution Strategy.

The best discrete-recombination certificate found in the recorded run uses again
\[
R=\frac{6672416}{100000}=66.72416,
\mbox{ with}\]
\[
\begin{split}
S_Q=\{&2,3,5,47,71,79,97,101,107,109,139,151,163,167,179,\\
      &191,211,223,239,241,251,263\}.
\end{split}
\]
The multiplicities are
\[
\begin{array}{c|rrrrrrrrrrr}
p&2&3&5&47&71&79&97&101&107&109&139\\
\hline
k(p)&49&29&19&8&7&7&6&6&6&6&6
\end{array}
\]
\[
\begin{array}{c|rrrrrrrrrrr}
p&151&163&167&179&191&211&223&239&241&251&263\\
\hline
k(p)&6&5&6&5&5&5&5&5&5&5&5
\end{array}
\]
\begin{figure}[htbp]
\centering
\resizebox{0.98\linewidth}{!}{%
\begin{tikzpicture}[
  font=\scriptsize,
  tcell/.style={draw=blue!70!black, fill=blue!9, rounded corners=1.5pt, minimum width=0.70cm, minimum height=0.47cm, align=center},
  inertcell/.style={draw=green!45!black, fill=green!10, rounded corners=1.5pt, minimum width=0.90cm, minimum height=0.70cm, align=center},
  ramcell/.style={draw=orange!75!black, fill=orange!18, rounded corners=1.5pt, minimum width=0.90cm, minimum height=0.70cm, align=center},
  basecell/.style={draw=gray!70!black, fill=gray!15, rounded corners=1.5pt, minimum width=0.35cm, minimum height=0.25cm, align=center},
  note/.style={align=left, font=\scriptsize},
  lab/.style={font=\scriptsize\bfseries, anchor=east}
]
\node[font=\small\bfseries, align=center] at (5.55,1.05) {Best discrete-recombination certificate};
\node[note, align=center] at (5.55,0.55) {$R=6672416/100000$, $\delta=0.0152628688\ldots$};

\node[lab] at (-0.25,0) {$T$};
\foreach \p [count=\i from 0] in {3,5,7,11,13,17,19,23,29,31,37,41,43} {
  \node[tcell] at (\i*0.88,0) {\p};
}
\node[note, anchor=west] at (0,-0.55) {The same 13-prime set $T$ is used; exactly seven entries are $3\bmod 4$.};

\node[lab] at (-0.25,-1.35) {$S_Q$};
\foreach \p/\k/\sty [count=\i from 0] in {2/49/ramcell,3/29/ramcell,5/19/ramcell,47/8/inertcell,71/7/inertcell,79/7/inertcell,97/6/inertcell,101/6/inertcell,107/6/inertcell,109/6/inertcell,139/6/inertcell} {
  \node[\sty] at (\i*1.02,-1.35) {\textbf{\p}\\$k=\k$};
}
\foreach \p/\k/\sty [count=\i from 0] in {151/6/inertcell,163/5/inertcell,167/6/inertcell,179/5/inertcell,191/5/inertcell,211/5/inertcell,223/5/inertcell,239/5/inertcell,241/5/inertcell,251/5/inertcell,263/5/inertcell} {
  \node[\sty] at (\i*1.02,-2.32) {\textbf{\p}\\$k=\k$};
}
\node[note, anchor=west] at (0,-2.95) {22 selected primes; orange means ramified in $Q$, green means inert in $Q$; no selected prime splits in $Q$.};

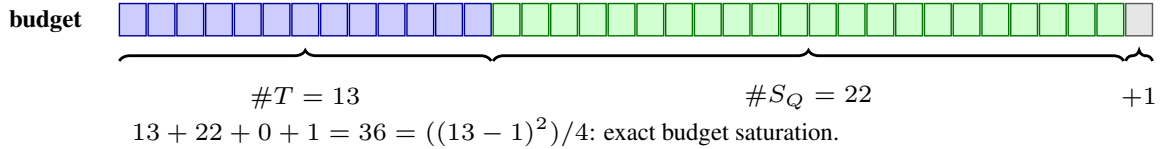
\node[lab] at (-0.25,-3.85) {budget};
\foreach \i in {0,...,12} {\draw[draw=blue!60!black, fill=blue!20] ({0.285*\i},-4.00) rectangle ++(0.265,0.32);}
\foreach \i in {13,...,34} {\draw[draw=green!45!black, fill=green!18] ({0.285*\i},-4.00) rectangle ++(0.265,0.32);}
\draw[draw=gray!70!black, fill=gray!20] ({0.285*35},-4.00) rectangle ++(0.265,0.32);
\draw[decorate, decoration={brace, amplitude=3pt}, thick] (0,-4.25) -- ({0.285*13-0.02},-4.25) node[midway, yshift=-0.36cm] {$\#T=13$};
\draw[decorate, decoration={brace, amplitude=3pt}, thick] ({0.285*13},-4.25) -- ({0.285*35-0.02},-4.25) node[midway, yshift=-0.36cm] {$\#S_Q=22$};
\draw[decorate, decoration={brace, amplitude=3pt}, thick] ({0.285*35},-4.25) -- ({0.285*36-0.02},-4.25) node[midway, yshift=-0.36cm] {$+1$};
\node[note, anchor=west] at (0,-4.95) {$13+22+0+1=36=((13-1)^2)/4$: exact budget saturation.};

\node[ramcell, minimum width=0.35cm, minimum height=0.25cm] at (0,-5.65) {};
\node[anchor=west] at (0.27,-5.65) {ramified in $Q$};
\node[inertcell, minimum width=0.35cm, minimum height=0.25cm] at (2.15,-5.65) {};
\node[anchor=west] at (2.42,-5.65) {inert in $Q$};
\node[tcell, minimum width=0.35cm, minimum height=0.25cm] at (4.05,-5.65) {};
\node[anchor=west] at (4.32,-5.65) {prime in $T$};
\node[basecell] at (5.80,-5.65) {};
\node[anchor=west] at (6.07,-5.65) {constant $+1$ budget term};
\end{tikzpicture}}
\caption{Parameter-level visualization of the best certificate found by the Tailored Integer Evolution Strategy with discrete recombination.  The figure records the finite data that are checked by the verification script; it is not a coordinate realization of the corresponding asymptotic planar construction.}
\label{fig:recomb-candidate-dashboard}
\end{figure}

The arithmetic checks again pass.  Formula~\eqref{eq:delta-objective} gives
\[
\begin{aligned}
\text{numerator}   &=4.5232596663564752\ldots,\\
\text{denominator} &=296.35710825977047\ldots,\\
\delta             &=0.01526286881700719\ldots .
\end{aligned}
\]
This is a small improvement over the non-recombining Tailored Integer Evolution Strategy certificate, whose value was $0.01526166106841929\ldots$.

\begin{remark}[Algorithmic side experiments]
Discrete recombination was beneficial in the recorded run, but only marginally.  A self-adaptive mutation-rate variant was also tested; it did not produce a better certificate than the simpler Tailored Integer Evolution Strategy, apparently because the strategy parameter converged prematurely to a narrow local optimization regime.  For this reason, the self-adaptive variant is not used as the main reported optimization result.
\end{remark}

\begin{remark}[Status of the claim]
The verification attempts passed all checks implemented in the self-contained implementation: primality, parity of $T$, admissibility of the primes in $S_Q$, absence of split primes in $Q$, exact saturation of the Golod--Shafarevich budget, positivity of all $k(p)$, $R>1$, and the numerical inequalities for the stated clean exponents.  Thus the bundle supports a sequence of certificate-level lower-bound improvements relative to Sawin's published $n^{1.014}$ statement: first the greedy certificate, then the Tailored Integer Evolution Strategy certificate, and finally the discrete-recombination variant.  However, the sharper displayed decimal values $1.0151718056\ldots$, $1.0152616611\ldots$, and $1.0152628688\ldots$ should be treated as candidate decimals until independently checked with interval arithmetic and reviewed by a human expert in the number-theoretic construction.
\end{remark}
\newpage
\section*{Dedication}

\textit{Dedicated to the memory of G\"unter Rudolph}.

\section*{Acknowledgements}

The author thanks the researchers and engineers who created and openly released the 2026 counterexample material, and Noga Alon, Thomas F. Bloom, W. T. Gowers, Daniel Litt, Will Sawin, Arul Shankar, Jacob Tsimerman, Victor Wang, and Melanie Matchett Wood for the human-verified exposition.  The author is especially grateful to Will Sawin for the explicit quantitative lower bound and for making the associated parameter-selection problem visible.  Any remaining errors or overinterpretations are the author's responsibility.

\section*{Declaration on the Use of Artificial Intelligence}

OpenAI ChatGPT 5.5 was used as an auxiliary tool for programming assistance, code review, debugging, and cross-checking the interpretation and implementation of the constraints and optimization model in Sawin~\cite{Sawin2026}. In particular, AI-assisted dialogue supported understanding aspects of the number theoretic constraint system in Sawin~\cite{Sawin2026}.

The tailored integer evolution strategy, the greedy construction algorithm, the mathematical formulation, the optimization methodology, the computational experiments, the verification procedures, the related-work discussion, the selection and interpretation of references, and all reported certificates were designed, written, executed, and independently validated by the author and collaborators. Artificial-intelligence tools were not used to design the optimization algorithms, generate mathematical proofs, generate the certificates, prepare the related work discussion, or select references. Responsibility for the correctness of the implementation, computations, interpretations, references, and conclusions remains entirely with the author.

\end{document}